\documentclass[letterpaper, 12pt,onecolumn]{article}

\usepackage[utf8]{inputenc}
\usepackage{amsfonts}
\usepackage{amssymb}
\usepackage{mathrsfs}
\usepackage{amsmath}
\usepackage{amsthm}
		    
\def\cT{\mathcal T}

\def\tA{\tilde A}
\def\tr{^{\mathsf {T}}}
\def\eop{\unskip\nobreak\hfil\penalty50\hskip2em\hbox{}\nobreak
\hfill\mbox{$\Box $}\par}
\newtheorem{thm}{Theorem}[section]
  
\begin{document}
	  \title{Elementary triangular matrices and inverses of $k$-Hessenberg and triangular matrices}
    \author{Luis Verde-Star}
\date{Oct. 2, 2015}
\maketitle

	        \begin{abstract}
			{We use elementary  triangular matrices to obtain some factorization, multiplication, and inversion properties of triangular matrices. We also obtain explicit expressions for the inverses of  strict $k$-Hessenberg matrices and banded matrices. Our results can be extended to the cases of block triangular and block Hessenberg matrices.   An $n \times n$  lower triangular matrix is called elementary if it is of the form $I +C$,   where $I$ is the identity matrix and $C$ is  lower triangular and  has all of its nonzero entries in the $k$-th column, where $1 \le k \le n$. 	
			}
		\end{abstract}

\medskip
		  {\bf Keywords:} Triangular matrices, factorization, $k$-Hessenberg matrices, matrix inversion.

\section{Introduction}
The importance of triangular, Hessenberg, and banded matrices is well-known. Many  
 problems  in linear algebra and matrix theory are solved by some kind of reduction to problems involving such types of matrices. This occurs, for example, with the $LU$ factorizations and the $QR$ algorithms. 

In this paper, we study first some simple properties of triangular matrices using a particular class of such matrices that we call elementary. Using elementary matrices we obtain factorization and inversion properties  and a formula for powers of triangular matrices. 

Some of our results may be useful to develop parallel algorithms to compute powers and inverses of triangular matrices and also of block-triangular matrices.

In the second part of the paper we obtain an explicit formula for the inverse of a  strict  $k$-Hessenberg matrix in terms of the inverse of an associated triangular matrix. 
Our formula is obtained by extending $n \times n$ $k$-Hessenberg matrices to $(n+k) \times (n+k)$ invertible triangular matrices and   using  some natural block decompositions. Our formula  can be applied to find the inverses of tridiagonal and banded matrices.

The problem of finding the inverses of Hessenberg and banded matrices has been studied by several authors, using different approaches, such as determinants and recursive algorithms. See \cite{Elouafi}, \cite{Ikebe}, \cite{Piff}, \cite{Xu}, and \cite{Yama}.
  
\section{Elementary triangular matrices}
In this section $n$ denotes a fixed positive integer,  $N=\{1,2,\ldots,n\}$, and $\cT$ denotes the set of lower triangular $n \times n$ matrices with complex entries.
An element of $\cT$  is called {\em elementary} if it is of the form $I +C_k$,  for some $k \le n $,  where $I$ is the identity matrix and $C_k$ is  lower triangular   and  has all of its nonzero entries in the $k$-th column.

Let $A=[a_{j,k}]  \in \cT$.  For each $ k \in N$  we define $E_k$ as the matrix obtained from the  identity matrix $I_n$ by replacing its $k$-th column with  the $k$-th  column of $A$, that is, 
$(E_k)_{j,k}=a_{j,k}$ for $j=k,k+1,\ldots, n$, $(E_k)_{j,j}=1 $ for $ j \ne k$, and all the other entries of $E_k$ are zero. 
The matrices $E_k$ are called the {\em elementary factors} of $A$ because
 $$ A = E_1 E_2 \cdots E_n . \eqno(2.1)$$
 Let us note that performing the multiplications in (2.1) does not require any arithmetical operations. It is just putting the columns of $A$ in their proper places. 

 If for some $k$ we have $a_{k,k}\ne 0$ then $E_k$ is invertible and it is easy to verify that 
 $$  (E_k)^{-1} = I - \frac{1}{a_{k,k}} (E_k -I). \eqno(2.2)$$
 Note that $(E_k)^{-1}$ is also an elementary lower triangular matrix. 
 If $A$ is invertible then all of its elementary factors are invertible and from (2.1) we obtain
  $$  A^{-1} = (E_n)^{-1} (E_{n-1})^{-1} \cdots (E_2)^{-1} (E_1)^{-1}. \eqno(2.3) $$
  Therefore  $A^{-1}$  is the product of the elementary factors of the  matrix  $B= (E_1)^{-1} (E_2)^{-1}\cdots (E_n)^{-1}$, but in reverse order.

  Notice that $B= I + (I -A) D^{-1},$  where $D=\hbox{\, Diag}(a_{1,1},a_{2,2},\ldots,a_{n,n})$.
 Therefore
$$B^{-1}=E_n E_{n-1} \cdots E_2 E_1. \eqno(2.4)$$

  If we are only interested in computing the inverse of $A$ we can find the inverse of $\tA=A D^{-1}, $  which has all its diagonal entries equal to one, and then we have  $A^{-1}= D^{-1} {\tA}^{-1}$. This means that it would be enough to consider matrices with diagonal entries equal to one to study the construction of inverses of triangular matrices. If we consider general triangular matrices we can also obtain results about the computation of powers of such matrices. 
  We will obtain  next some results about products of elementary factors.

  We define $C_k= E_k - I$  for $k \in N$. Note that  $(C_k)_{k,k}= a_{k,k} -1$,  $(C_k)_{j,k}= a_{j,k} $ for $j=k+1,\ldots,n$, and all the other entries are zero, that is, all the nonzero entries of $C_k$ lie on the $k$-th column.
 It is clear that $A$ has the following additive decomposition
 $$A= I + C_1+C_2+ \cdots + C_n. \eqno(2.5)$$

 Let $L$ be the $n \times n$ matrix such that $L_{k+1,k}=1$ for $k=1,2,\ldots,n-1$ and all its other entries are zero. We call $L$ the  shift matrix,  since  $M L$ is $M$ shifted one column to the left.
Note that $\cT$ is invariant under the maps $M \to M L$   and $M \to L M$.
  
 In order to simplify the notation we will write the entries in the main diagonal of $A$ as  $a_{k,k}=x_k$ for $k=1,2,\ldots,n.$
 
The matrices $C_k$ have simple multiplication properties that we list in the following proposition.  
 \begin{thm}

	 \begin{enumerate}
		 \item{ If $ j < k $ then $C_j C_k=0$.}
		 \item{ For $m \ge 1$ we have $C_k^m=(x_k -1)^{m-1} C_k$.}
		 \item{ If $ k > j  $  then $ C_k C_j= a_{k,j} C_k L^{k-j}.$}
	 \end{enumerate}
 \end{thm}
 The proofs are straightforward computations.
 Notice that all the nonzero entries of $C_k L^{k-j}$ are in the  $j$-th column.

From part 3 of the previous theorem we obtain immediately the following multiplication formula. If  $r\ge 2$  and  $k_1 >k_2 > \cdots > k_r$ then  
$$ C_{k_1} C_{k_2} \cdots C_{k_r}=a_{k_1,k_2} a_{k_2,k_3} \cdots a_{k_{r-1},k_r} C_{k_1} L^{k_1 -k_r}. \eqno(2.6)$$

If $K$ is a  subset of $N$   with at least two elements we define the matrix
$$ G(K)= a_{k_1,k_2} a_{k_2,k_3} \cdots a_{k_{r-1},k_r} C_{k_1} L^{k_1 -k_r}, \eqno(2.7)$$
where $K=\{k_1,k_2,\ldots, k_r\}$ and $k_1 >k_2 > \cdots > k_r.$ 
 Let us note  that all the nonzero entries of $G(K) $  are in its  $k_r$-th column.
If $K$ contains only one element, that is, $K=\{ j \}$  then we put 
 $G(K)=C_j$.

 Since $ E_k= I + C_k$, the multiplication properties of the $C_k$ can be used to obtain some corresponding  properties of the $E_k$. 
 
 \begin{thm}

	          \begin{enumerate}
			  \item{ If $ k > j  $ then $E_k  E_j = I + C_k +C_j + G(\{k,j\}).$}
			 \item{ For $m \ge 1$ we have $$E_k^m=I  + \frac{ x_k^m -1}{x_k -1}\,   C_k.   $$}
			 \item{ If   $K=\{k_1,k_2,\ldots, k_r\}$ and $k_1 >k_2 > \cdots > k_r$   then 
     $$ E_{k_1} E_{k_2} \cdots E_{k_r} = I + \sum_{J \subseteq  K} G(J). $$ }  
	 \end{enumerate}
       \end{thm}
 Proof: The proof of the first part is trivial. For the second part, use part 2 of Theorem 2.1 and the binomial formula. For part 3, write each factor in the form  $E_{k_j}=I + C_{k_j}$, expand the product, collect terms and use the definition of the function $G$.
Alternatively, we can  use part 1 repeatedly and then use the definition of $G$.
 \eop

Observe that the number of summands in part 3 is at most equal to the number of subsets of $K$. If some of the entries $a_{j,k}$ are equal to zero then some of the matrices $G(J)$ are zero.

Taking $K=N$ in part 3 of Theorem 2.2 we obtain
$$ E_n E_{n-1} \cdots E_2 E_1 = I + \sum_{J\subseteq N} G(J) . \eqno(2.8) $$
Let $k \in N $. Then the summands in (2.8) that may have nonzero entries in the $k$-th column (other than $I$) are the matrices $G(\{k\}\cup J)$ where $J \subseteq \{k+1,k+2,\ldots,n\}$, and thus the number of such matrices is at most equal to $2^{n-k}$. For $k=n$ there is only one, which is $C_n$, for $k=n-1$ they are $C_{n-1}$ and $G(\{n-1,n\})$, and so on. 
Since $G(\{k\}\cup J)$ is a scalar multiple of $ C_m L^{k-m} $ where $m$ is the largest element of $J$, we can group together all the terms having the same largest element $m$ and therefore the $k$-th column of the product in (2.8) is a linear combination of $C_k, C_{k+1}L, C_{k+2}L^2, \ldots, C_n L^{n-k}$, and the $k$-th column of the identity matrix.

Therefore, if we have computed $E_n E_{n-1} \cdots E_{k}$ then the columns with indices $k,k+1,\ldots,n$ are determined and do not change when we multiply by $E_{k-1}, E_{k-2}$, etc. Thus  $E_n E_{n-1} \cdots E_{k+1}$  and $E_n E_{n-1} \cdots E_{k+1} E_{k}$ only differ in their $k$-th columns. This means that computing the sequence $E_n E_{n-1} \cdots E_{k}$, for $k=n, n-1,\ldots,2,1$ we obtain $B^{-1}$ column by column, starting from the last one. 
This procedure may be useful to develop parallel algorithms for the computation of inverses of triangular matrices.

We consider now explicit expressions for the positive powers of a triangular matrix $A$ in terms of its elementary factors. We start with $A^2$.  Using (2.5) and part 1 of Theorem 2.1  we get
$$A^2=(I + C_1+C_2+\cdots + C_n)^2=I+ 2 \sum_{k=1}^n C_k +\sum_{k=1}^n C_k^2 + \sum_{j<k} C_k C_j. \eqno(2.9)$$
Therefore, by Theorem 2.1 we have 
$$A^2=I + \sum_{k=1}^n (1+ x_k) C_k +\sum_{j<k} a_{k,j} C_k L^{k-j} , \eqno(2.10)$$
where the last sum runs over all pairs $(j,k)$ such that $1 \le j < k \le n$.

Let $K=\{k_1,k_2,\ldots, k_r\}\subseteq N$, where  $k_1 >k_2 > \cdots > k_r, $ and let $m$ be a positive integer.
We define the scalar valued function $g(K,m)$ as follows
$$g(K,m)=\Delta[1, x_{k_1},x_{k_2},\ldots, x_{k_r}] t^{m+r}, \eqno(2.11)$$
where $\Delta[1, x_{k_1},x_{k_2},\ldots, x_{k_r}]$ denotes the divided differences functional with respect to the numbers $1, x_{k_1},x_{k_2},\ldots, x_{k_r}$. 
 $g(K,m)$ is a symmetric polynomial in the $x_j$. For the properties of divided differences see \cite{ddlrs}.

Using induction and some basic properties of divided differences we can  obtain an expression for $A^m$ in terms of matrices of the form $G(J)$.
\begin{thm}  For $m \ge 1$ we have 
$$A^m= I + \sum_{j=1}^n g(\{j\},m-1) C_j +\sum_{J\subset N,\ |J|=2} g(J,m-2) G(J)+ \qquad \qquad \qquad $$
$$\qquad  \qquad \qquad  \sum_{J\subset N,\ |J|=3} g(J,m-3) G(J)+ \cdots + \sum_{J\subset N,\ |J|=m} g(J,0) G(J). \eqno(2.12)$$
\end{thm}

Note that  the numbers  $g(J,0)$  in the last sum are all equal to 1.
A similar result for triangular matrices with distinct diagonal entries $x_j$ was obtained by Shur \cite{Shur}.
Other formulas for powers of general square matrices appear in \cite{FoM}. 
 
\section{Inverses of Hessenberg matrices}
In this section we use the fact that a  $k$-Hessenberg matrix $H$  is a submatrix of a larger triangular matrix to obtain a formula for the inverse of $H$, in case it exists.
We also characterize the invertible $k$-Hessenberg matrices in terms of properties of the  blocks in a natural block decomposition.

We call an   $n \times n$ matrix $H=[h_{i,j}]$ lower $k$-Hessenberg if  $h_{i,j}=0$ for $i < j-k$. We say that $H$ is strict lower $k$-Hessenberg if $h_{i,j}\ne 0$ for $i=j-k$. 

Any  lower $k$-Hessenberg  $n \times n$  matrix $H$ has the following block decomposition.
Let $m=n-k$. Then 
$$H= \left[ \begin{matrix} B & A \cr  D & C \cr \end{matrix} \right], \eqno(3.1)$$
where $A$ is $m \times m$, $B$  is $m \times k$, $C$ is $k \times m$, and $D$ is $k \times k$. Note that $A$ is lower triangular and it is invertible if $H$ is strict $k$-Hessenberg. 

Extending the block decomposition of equation (3.1) to form a lower triangular matrix we obtain the following result.

\begin{thm}
	Let $H$ be strict $k$-Hessenberg with the block decomposition given in (3.1).  Then $H$ is invertible if and only if $C A^{-1} B -D$ is invertible and if $H$ is invertible we have 
	$$ H^{-1} = \left[ \begin{matrix} 0 & 0 \cr A^{-1} & 0 \cr \end{matrix} \right] - \left[ \begin{matrix}   I_k \cr E \cr \end{matrix} \right]  G^{-1} \left[ \begin{matrix} F & I_k \cr \end{matrix} \right], \eqno(3.2) $$
		where $I_k$ is the $k \times k$ identity matrix,  $E=-A^{-1}B,$ \ \ $F=-C A^{-1}$,  and $G=C A^{-1} B -D$. 
	\end{thm}
Proof:  Suppose that $G=C A^{-1} B -D$ is invertible. Define the lower triangular $(n+k) \times (n+k)$ matrix  $T$  by
$$ T=  \left[ \begin{matrix}  I_k & 0 & 0 \cr B & A & 0 \cr D & C & I_k \cr \end{matrix} \right]. \eqno(3.3) $$ 
 It is easy to verify that 
 $$ T^{-1} =  \left[ \begin{matrix}  I_k & 0 & 0 \cr  E & A^{-1} & 0 \cr G & F & I_k \cr \end{matrix} \right], \eqno(3.4)$$ 
where $E, F, $ and $G$ are as previously defined. Consider now the block decomposition 
$$ T= \left[ \begin{matrix}  R & 0 \cr H & S \cr \end{matrix} \right], $$
	where $$ R=  \left[ \begin{matrix} I_k & 0 \end{matrix} \right], \qquad S= \left[ \begin{matrix}  0 \cr I_k \end{matrix} \right]. $$
		From  $T T^{-1}  = I_{n+k}$  we obtain the equations
		$$ H \left[ \begin{matrix} 0 & 0 \cr A^{-1} & 0 \cr \end{matrix} \right]   + S \left[ \begin{matrix} F & I_k \cr \end{matrix} \right] = I_n, \eqno(3.5)$$ 
	and
	$$ H \left[ \begin{matrix} I_k \cr E \cr \end{matrix} \right] +S G = 0. \eqno(3.6)$$
	Since $G$ is invertible by hypothesis,  we can solve for $S$ in the last equation and substitute the resulting expression for $S$  in equation (3.5). In this way we obtain 
	$$ H \left[ \begin{matrix} 0 & 0 \cr A^{-1} & 0 \cr \end{matrix} \right]  - H \left[ \begin{matrix} I_k \cr E \cr \end{matrix} \right] G^{-1} \left[ \begin{matrix} F & I_k \cr \end{matrix} \right] = I_n , $$ 
	which implies that $H$ is invertible and   (3.2) holds.

 Now suppose that $H$ is invertible and let 
  $$ H^{-1}= \left[ \begin{matrix} U & V \cr W & Y \cr \end{matrix} \right] , $$
 where the block decomposition is compatible with the decomposition of $H$ given in (3.1).
 Then, from the equation $H^{-1} H =I_n$ we obtain  
 $$ UB +VD =I_k, \qquad  UA + VC =0. \eqno(3.7)$$
 Since $A$ is invertible the second  equation in (3.7) yields $U + VC A^{-1}=0$, and multiplying by $B$ on the right we obtain $ UB + V C A^{-1}B =0$. Combining this last equation with the first one in (3.7) we get  $ V ( C A^{-1}B -D) = -I_k$ and therefore $G$ is invertible and $G^{-1}=-V$. \eop

	Let us observe that an important part of the computation  of $H^{-1}$  is the computation of the inverse of the triangular matrix $A$.   
Note also that  any given strict $k$-Hessenberg matrix can be modified to become invertible by changing the $k \times k$ block $D$ of (3.1) in a suitable way.

In the case of  $k=1$ the matrix  $G$ reduces to a number and then the second term in the right-hand side of (3.2) is the product of a column vector times a row vector. See \cite{Ikebe}.

 Note that Theorem 3.1 holds for tridiagonal matrices and also for banded matrices. 
 Suppose that $k=1$, $H$ is tridiagonal, and $n \ge 3$. Then  the matrices in the block decomposition of $H$ 
 are $B=[ h_{1,1}\ h_{2,1} \ 0 \ 0 \ \cdots \ 0]\tr$, \  $C=[0 \ 0 \ \cdots \ 0\ h_{n,n-1}\ h_{n,n}],$  and $D=0$. In this case $A$ is lower triangular and tridiagonal. Using the row version of the theory of elementary triangular matrices, that we describe in the next section, it is easy to construct a recursive algorithm to compute the inverses of tridiagonal matrices as the size $n$ increases.

The proof of Theorem 3.1 only uses the hypothesis that the block $A$ is invertible. Therefore the theorem holds also for other types of matrices, such as block Hessenberg matrices.

	  \section{Row and  block versions of the theory 
 of elementary triangular matrices}
 In this section we present a brief description of two variations on the theory of elementary triangular matrices presented in section 2. The first one is obtained when we consider lower triangular matrices that are the sum of the identity matrix plus a matrix that has  all of its nonzero entries in a single {\em row.} In the second one we consider block lower triangular matrices with square blocks along the diagonal that may have different sizes.

 Recall that $\cT$ denotes the $n \times n$ lower triangular matrices.  
 An element of $\cT$  is called {\em row elementary} if it is of the form $I +R_k$,  for some $k \le n $,  where $I$ is the identity matrix and $R_k$ is  lower triangular   and  has all of its nonzero entries in the $k$-th row. 

Let $A=[a_{k,j}]  \in \cT$.  For $ k \in N$  we define $F_k$ as the matrix obtained from the  identity matrix $I$ by replacing its $k$-th row  with  the $k$-th  row  of $A$, that is, 
$(F_k)_{k, j}=a_{k, j}$ for $j=1,2,\ldots, k$, $(F_k)_{j,j}=1 $ for $ j \ne k$, and all the other entries of $F_k$ are zero.
The matrices $F_k$ are called the {\em elementary factors by rows } of $A$ because
 $$ A = F_1 F_2 \cdots F_n . \eqno(4.1)$$
If for some $k$ we have $a_{k,k}\ne 0$ then $F_k$ is invertible and 
 $$  (F_k)^{-1} = I - \frac{1}{a_{k,k}} (F_k -I). \eqno(4.2)$$
 Therefore, if $A$ is invertible then
 $$ A^{-1}= (F_n)^{-1} (F_{n-1})^{-1} \cdots (F_2)^{-1} (F_1)^{-1}. \eqno(4.3)$$
The matrices $(F_j)^{-1}$ are the elementary factors by rows of the matrix 
$$B=  I - \frac{1}{a_{k,k}} (A -I), $$  and 
$$B^{-1}=F_n F_{n-1} \cdots F_2 F_1. \eqno(4.4)$$

 Define $ R_k =F_k -I$  for $k \in N$. Then $A=I + R_1+R_2+\cdots +R_n$. 
It is easy to see that $F_k F_{k-1} \cdots F_2 F_1$ and $F_{k-1} \cdots F_2 F_1$ only differ in the $k$-th row, and the difference is a linear combination of translates of $R_1, R_2, \ldots, R_k$. Note that $F_k F_{k-1} \cdots F_2 F_1$ is the inverse of the submatrix of $B$ obtained by deleting the rows and columns with indices $k+1,k+2, \ldots, n$, which is often called the $k \times k$ section of $B$. Therefore, computing the sequence of matrices $F_k F_{k-1} \cdots F_2 F_1$ for $k=1,2,3,\ldots $  yields a recursive algorithm that gives the inverse of $B$ row by row. That algorithm can also be used to find the inverses of the sections of infinite lower triangular matrices such as the ones considered in \cite{LVS}. The inversion algorithm introduced in \cite{LVS} can be combined with 
the computation of inverses of diagonal blocks of a triangular matrix, using multiplication of elementary matrices, by rows or by columns. 

The concept of elementary triangular matrices (by colummns of rows) can be generalized to the case of block triangular matrices in a natural way. We describe next how it is done in the case of column block elementary matrices.

Let $k_1,k_2,\ldots,k_r$ be positive integers such that $n=k_1+k_2+\cdots +k_r$.  Let $X_j$ be a $k_j \times k_j$ matrix for $1 \le j \le r$ and let $A$ be an $n \times n$ block matrix that has the matrices $X_j$  along the diagonal and all its other nonzero entries below the diagonal blocks.  

For $j \in \{1,2,\ldots,r\}$ let $E_j$, called the block elementary factor of $A$  
by columns, be the matrix that coincides with $A$ in all the columns corresponding to the diagonal block $X_j$, that is, the columns with indices between $k_1+k_2+\cdots +k_{j-1}+1$  and $k_1+k_2+\cdots +k_{j}$,  and coincides with the identity matrix in the rest of the columns. Then we have   $A= E_1 E_2 \cdots E_r$. 

If the block $X_j$ is invertible then $E_j$ is also invertible and 
$$(E_j)^{-1}= I - (E_j-I)  \hbox{\,Diag}(I_{m_j}, (X_j)^{-1}, I_{p_j}),\eqno(4.5)$$
where $ \hbox{\,Diag}(I_{m_j}, (X_j)^{-1}, I_{p_j})$ is the block diagonal matrix that coincides with the block diagonal of $A$ in the $j$-th block and with the identity matrix  in the rest of the blocks. Note that $m_j=k_1+k_2+ \cdots+ k_{j-1}$ and $p_j = n- m_j - k_j$. 
If all the diagonal blocks $X_j$ are invertible then 
$$A^{-1}= (E_r)^{-1} (E_{r-1})^{-1}\cdots (E_2)^{-1} (E_1)^{-1}. \eqno(4.6)$$

       Department of Mathematics, Universidad Aut\'onoma Metropolitana,\\
 Iztapalapa, Apartado 55--534, Mexico D.F., Mexico,\\ 
 e-mail: verde@xanum.uam.mx. 
\end{document}